\documentclass[a4paper,11pt,leqno]{article}
\usepackage[dvips]{graphicx}
    \usepackage[usenames]{color}

\usepackage{amsmath}
\usepackage{amsfonts}
\usepackage{amssymb}

\newtheorem{tw}{Theorem}

\newtheorem{lem}[tw]{Lemma}

{

\newcommand{\xx}[1]{\bar{x_{#1}}}
\begin{document}
\def\pf{{\bf Proof. }}
\def\thetw{\thesection.\arabic{tw}.}

\begin{center}
{\LARGE Viscosity iteration in CAT($\kappa$) spaces}

\vspace{5mm}

\fontencoding{T1}

\Large{Bo\.zena PI\k{A}TEK}\\

\fontencoding{OT1}
\vspace{2mm}

\small{Institute of Mathematics, Silesian University of Technology, 44-100 Gliwice, Poland,
{\it b.piatek@polsl.pl}}

\end{center}

\vspace{5mm}

{\small
\noindent{\bf Abstract} We study the approximation of fixed points of nonexpansive mappings in CAT($\kappa$) spaces. We show that the iterative sequence generated by the Moudafi's viscosity type algorithm converges to one of the fixed points of the nonexpansive mapping depending on the contraction applied in the algorithm.
\\}


\pagebreak

\section{Introduction}

The problem of existence of fixed points for different types of mappings is one of the most important in the mathematical analysis. If a~mapping $T\colon X\to X$, where $X$ is a~complete metric space, is a~contraction then the well known Banach contraction principle implies that there is a unique fixed point of $T$, i.e., a~point $x_0$ for which $T(x_0)=x_0$, and moreover each Picard iterative sequence $(T^nx)$ tends to $x_0$. Much more complicated is the case when we only know that $T$ is nonexpansive, i.e.,
$d(Tx,Ty) \leq d(x,y)$ for all couples of points of $X$. There is a~huge literature on the theory of fixed points for nonexpansive mappings, for comprehensive expositions the reader can check \cite{goebelkirk, handbook} and references therein.

In the case that at least one fixed point exists, then a very important problem is to know how to locate or approximate the fixed point. In 2000 Moudafi \cite{moudafi} proposed one of the most successful methods in which a~certain iterative sequence strongly convergences to fixed points. More precisely, let $f\colon X\to X$ be a~contraction. The iterative procedure is done by the~formula:
$$
x_{n+1}=t_n f(x_n)+(1-t_n) T(x_n)
$$
for a~fixed sequence $(t_n)$ of real numbers. Then, under suitable conditions including that $X$ is a Hilbert space, the above iterative sequence tends to the fixed point for which the distance to the image of $Fix T$ under $f$ is minimal. Clearly, by $Fix T$ we understand the set of fixed points of $T$.
A~few years later this viscosity method was modified by Nilsrakoo and Saejung in \cite{saejung} who considered results of two--steps approximation due to Chidume and Chidume (see \cite{chidume}). This modification allows to weaken the assumptions on the sequence of real numbers $(t_n)$. However, the problem has only been considered so far for Hilbert or uniformly convex Banach spaces. In this work we take up the problem of studying this two-steps viscosity scheme in CAT($\kappa$) spaces. This class of spaces form a~natural and vast extension of spherical geometry but not only. CAT($\kappa$) spaces generalizes complete Riemannian manifolds with curvature bounded above, Tits buldings and limits of Hadamard spaces, among several others. We will particularized for $\kappa=1$ and then we will obtain results for any positive $\kappa$ and, by a limiting process, also for $\kappa=0$. Since our approach to the problem focuses on CAT($1$) spaces we will need to develop a number of technical results which cannot be directly inherited from previous works on Hilbert or Banach spaces. 

Previous results about approximation of fixed points for CAT($0$) spaces were first obtained in 2008 by Dhompongsa and Panyanak in \cite{Dhompongsa}. Then more results on approximation of fixed points for CAT($0$) spaces were considered by different authors (see e.g. \cite{Khamsi,Laowang,Halpern,Chou}). In \cite{Khan} we can find the most alike algorithm method developed in CAT($0$) spaces to the one we study here. Very recently, in \cite{moja} the author proposed the application of the Halpern iteration scheme to the more general case of CAT($\kappa$) spaces with $\kappa >0$, see also the still more recent work \cite{gena}.

The main goal of this paper is a~generalization of \cite{moja} by an~application of a viscosity type iterative method. More precisely, we prove that a similar iterative procedure as that one in \cite{saejung} gives us strong convergence to the fixed point under suitable conditions in CAT($\kappa$) spaces. To achieve our result we do not follow the usual approach in Banach but rather we attend to the geometric properties of CAT($\kappa$) spaces.

Our work is organized in the following way. In Section 2 we
introduce some definitions and notations that will be needed in the
rest of the work. In Section $3$ we present some lemmas devoted to the geometry of CAT($\kappa$) spaces, e.g. the modification of quadrilateral cosine formula for CAT($\kappa$) spaces proposed by Berg and Nikolaev in \cite{berg}. In Section 4 we prove the main results of the paper, i.e., the convergence of the viscosity iterative sequence in complete CAT($\kappa$). We first obtain it for $\kappa >0$ and then deduce the result for CAT($0$) spaces. Our iterative procedure is a much more complicated one than others already studied in this setting and so it will request from an extra condition on the geometry of CAT($\kappa$) spaces, that is, the N-property. This property was first introduced in \cite{RafaAurora} and, as far as the author knows, no example of a CAT($0$) space has been exhibited in the literature to miss property N yet. We close this work with an appendix where we show an example of such a space, that is, a CAT($0$) space lacking property N.

\section{Preliminaries}

\setcounter{equation}{0}\setcounter{tw}{0}

Let $(X,\rho)$ be a metric space. Suppose that for each pair of points $x,y\in X$ there is an isometric embedding $\varphi\colon[0,d(x,y)]\to X$, such that $\varphi (0)=x$, $\varphi (d(x,y))=y$. Then the image $\varphi ([0,d(x,y)])$ is called {\it a metric segment} which, when unique, is usually denoted by $[x,y]$. Such a~space $X$ is called a~geodesic space.
In the sequel we assume that $(X,\rho)$ is geodesic.
%Let $x,y,z\in X$. 
{\it A geodesic triangle} $\triangle(x,y,z)$ consists of triples of points $x,y,z$ (called vertices of the triangle) and three metric segments joining each pair of these vertices (called its edges). 

Now we introduce the concept of model spaces $M_\kappa^2$, for $\kappa \geq 0$, which we will
need to recall the definition of $\rm{CAT}(\kappa)$ spaces. Very generous expositions on this type of spaces may be found in \cite{bridson} and \cite{burago}. 

Let us denote by $\left<\cdot,\cdot\right>$ the scalar product in $\mathbb{R}^3$. By $\mathbb{S}^2$ we denote the unit sphere in $\mathbb{R}^3$, i.e., the set $\{(x_1,x_2,x_3)\in\mathbb{R}^3 : x_1^2+x_2^2+x_3^2=1\}$. In $\mathbb{S}^2$ one may introduce a metric (called the spherical distance, see \cite[section
3.1]{goebelreich}) in such a way that
$$
d_{\mathbb{S}^2}(x,y)=\arccos\left<x,y\right>.
$$
In the sequel this metric will be denoted by $d$ and $\rho$ will state for the metric in a general geodesic metric space.
Moreover, on the unit sphere the spherical cosines law holds
$$
\cos d(y,z)=\cos d(x,y)\cos d(x,z)+\sin d(x,y)\sin d(x,z)\cos \alpha,
$$
where $\alpha$ denotes the spherical angle between geodesic $[x,y]$ and $[x,z]$ on $\mathbb{S}^2$.

The model space $M_\kappa^2$ is defined as follows
\begin{enumerate}
	\item[(i)] the Euclidean plane $\mathbb{E}^2$ in case of $\kappa=0$;
	\item[(ii)] the unit sphere with distances multiplying by $\frac{1}{\sqrt{\kappa}}$ in case of $\kappa>0$.
\end{enumerate}

Let $\triangle(x_1,x_2,x_3)$ be a geodesic triangle in a geodesic space $(X,\rho)$. The triangle $\triangle(\bar{x_1},\bar{x_2},\bar{x_3})$ from $M_\kappa^2$ is called {\it a comparison triangle} of $\triangle(x_1,x_2,x_3)$ if $\rho(x_i,x_j)=d(\bar{x_i},\bar{x_j})$ states for all pairs $i,j\in\{1,2,3\}$, where $d$ denotes the metric in the model space $M_\kappa^2$. It is well known that given a~geodesic
triangle then a~comparison triangle exists if
$\rho(x_1,x_2)+\rho(x_2,x_3)+\rho(x_3,x_1)<2D_\kappa$, where
$D_\kappa=\frac{\pi}{\sqrt{\kappa}}$ for a~positive $\kappa$ and
$D_0=\infty$. Moreover, the comparison triangle is unique up to
isometry.

For each pair of points $y_1\in [x_1,x_2]$ and $y_2\in [x_1,x_3]$ one may find comparison points $\bar{y_1}\in [\bar{x_1},\bar{x_2}]$ and $\bar{y_2}\in [\bar{x_1},\bar{x_3}]$ in such a way that $\rho(x_1,y_i)=d(\bar{x_1},\bar{y_i})$, $i\in\{1,2\}$. 

The space $X$ is said to be a CAT($\kappa$) space if for each pair $(y_1,y_2)$ of points of any triangle $\triangle(x_1,x_2,x_3)\subset X$ for which $\rho(x_1,x_2)+\rho(x_2,x_3)+\rho(x_3,x_1)<2D_\kappa$ and for their comparison points $\bar{y_i}\in\triangle(\bar{x_1},\bar{x_2},\bar{x_3})$ {\it the CAT($\kappa$) inequality}
$$
\rho(y_1,y_2)\leq d(\bar{y_1},\bar{y_2}),
$$
holds.

Let $X$ be a~$\rm{CAT(\kappa)}$ space. By $\alpha x + (1-\alpha)y$
where $\alpha\in(0,1)$ and $x,y\in X$, $d(x,y)< D_\kappa$, we denote a~point of the
metric segment $[x,y]$ such that $\rho(x,\alpha
x+(1-\alpha)y)=(1-\alpha)\rho(x,y)$ and $\rho(\alpha
x+(1-\alpha)y,y)=\alpha \rho(x,y)$. Since that metric segment is unique (see \cite[Proposition 1.4 (1), p.~160]{bridson}) that
combination is well defined. Moreover, a~set $C\subset X$ is
called {\it $D_\kappa$-convex} when $[x,y]\subset C$ for all
$x,y\in C$ such that $\rho(x,y)<D_\kappa$. It is well-known that
nonempty closed and $D_\kappa$-convex subsets of complete
CAT($\kappa$) spaces are proximinal (e.g. \cite[Proposition 2.4,
Exercise 2.6(1), p. 176--178]{bridson}), that is, for $C$
a~nonempty closed and $D_\kappa$-convex subset of a CAT($\kappa$)
space $X$ and for any $x\in X$ such that
$\rho(x,C)<\frac{D_\kappa}{2}$ the set
$$
P_C(x):=\{ y\in C:\; d(x,y)=\inf_{z\in C}d(x,z)\}
$$
is a~singleton.

A CAT($\kappa$) space $X$ is said to satisfy {\it N--property} if for any closed $D_\kappa$-convex set $C\subset X$ and for each points $x_1,x_2\in X$ such that $\rho(x_i,C)<\frac{D_\kappa}{2}$, $\rho(x_1,x_2)<\frac{D_\kappa}{2}$ and $P_C(x_1)=P_C(x_2)=P$, the projection $P_C(\alpha x_1 + (1-\alpha)x_2)$, $\alpha\in(0,1)$, is also equal to $P$. As a consequence, if property N happens, then the continuity of the metric projection implies that $P_{[a,b]}(\alpha x_1 + (1-\alpha)x_2)\in [P_{[a,b]}(x_1),P_{[a,b]}(x_2)]$ for any metric segment $[a,b]$.

Let $X$ be a geodesic space. For each $\alpha,\beta\in(0,1)$ let us consider $\triangle(\bar{x},\bar{y_\alpha},\bar{z_\beta})$ the comparison triangle in the Euclidean plane of $\triangle(x,(1-\alpha)x+\alpha y,(1-\beta)x+\beta z)$. Then let $\angle_{\bar{x}}(\bar{y_\alpha},\bar{z_\beta})$ denotes the angle between edges $[\bar{x},\bar{y_\alpha}]$ and $[\bar{x},\bar{y_\alpha}]$, i.e.,
$$
\angle_{\bar{x}}(\bar{y_\alpha},\bar{z_\beta})=\arccos\dfrac{d^2(\bar{x},\bar{y_\alpha})+d^2(\bar{x},\bar{z_\beta})-d(\bar{y_\alpha},\bar{z_\beta})}{2d(\bar{x},\bar{y_\alpha})d(\bar{x},\bar{z_\beta})}.
$$
{\it The Alexandrov angle} $\angle_x(y,z)$ between metric segments $[x,y]$ and $[x,z]$ is defined as
$$
\angle_x(y,z)=\limsup\limits_{\alpha,\beta\to 0}\angle_{\bar{x}}(\bar{y_\alpha},\bar{z_\beta})
$$
(compare \cite[Definition I.1.12, p. 9]{bridson}). Let us emphasise that in case of $X$ being a~CAT($\kappa$) space,  the limit $\lim\limits_{\alpha,\beta\to 0}\angle_{\bar{x}}(\bar{y_\alpha},\bar{z_\beta})$ exists if only $\rho(x,y)+\rho(y,z)+\rho(z,x)<2D_\kappa$ (see \cite[Proposition II.3.1, p. 184]{bridson}).

\section{Technical results}

\setcounter{equation}{0}\setcounter{tw}{0}

In this section we present some technical results that will be needed in our main section.

%Before we show the main results of this paper some technical lemmas will be needed. The first three give us a technical estimation in case of CAT(1) spaces. The proofs of Lemmas 3.2 and 3.3 may be found in \cite{moja}.

\begin{lem}\label{lem:pierwszy}
 Let $x_1,x_2,x_3,x_4$ be different points of a~CAT(1) space $X$ such that $\rho(x_i,x_j)<\frac{\pi}{2}$, $i,j\in\{1,2,3,4\}$
and $x_3\in[x_2,x_4]$. Then
$$
\sin^2\frac{\rho(x_1,x_3)}{2} \leq
\frac{\rho(x_2,x_3)}{\rho(x_2,x_4)}\sin^2\frac{\rho(x_1,x_4)}{2}+\frac{\rho(x_3,x_4)}{\rho(x_2,x_4)}\sin^2\frac{\rho(x_1,x_2)}{2}.
$$
\end{lem}

\pf
Let $\Delta(\xx{1},\xx{2},\xx{4})$ be a comparison triangle of $\Delta(x_1,x_2,x_4)$ on a unit sphere. Consider the triangles $\Delta(\xx{1},\xx{2},\xx{3})$ and $\Delta(\xx{1},\xx{3},\xx{4})$ and apply the cosine law to each of them. By solving equal terms we obtain
$$
\frac{\cos d(\xx{1},\xx{3})}{\sin d(\xx{2},\xx{3})}-\frac{\cos d(\xx{1},\xx{2})\cos d(\xx{2},\xx{3})}{\sin d(\xx{2},\xx{3})}=\frac{\cos d(\xx{1},\xx{4})}{\sin d(\xx{2},\xx{4})}-\frac{\cos d(\xx{1},\xx{2})\cos d(\xx{2},\xx{4})}{\sin d(\xx{2},\xx{4})},
$$
what yields
$$
\cos d(\xx{1},\xx{3})=\cos d(\xx{1},\xx{4})\frac{\sin d(\xx{2},\xx{3})}{\sin d(\xx{2},\xx{4})}
$$
$$
+\cos d(\xx{1},\xx{2})\frac{\sin d(\xx{2},\xx{4})\cos d(\xx{2},\xx{3})-\sin d(\xx{2},\xx{3})\cos d(\xx{2},\xx{4})}{\sin d(\xx{2},\xx{4})}
$$
$$
=\cos d(\xx{1},\xx{4})\frac{\sin d(\xx{2},\xx{3})}{\sin d(\xx{2},\xx{4})}+\cos d(\xx{1},\xx{2})\frac{\sin d(\xx{3},\xx{4})}{\sin d(\xx{2},\xx{4})}
$$
$$
\geq\cos d(\xx{1},\xx{4})\frac{d(\xx{2},\xx{3})}{ d(\xx{2},\xx{4})}+\cos d(\xx{1},\xx{2})\frac{d(\xx{3},\xx{4})}{d(\xx{2},\xx{4})}.
$$
This, on account of the inequality $\rho(x_1,x_3)\leq d(\xx{1},\xx{3})$, completes the proof by recalling that $\cos 2x= 1-2\sin^2 x$ for trigonometric real functions.$\qquad\blacksquare$

\vskip.5cm

The next lemma is Lemma 3.3 in \cite{moja}.

\begin{lem}\label{CAT}
Let $\triangle(A,B,C)$ be a~triangle in a~CAT(1) space with all
sides no longer than $M$, where $M\leq \frac{\pi}{2}$. Moreover,
suppose that points $D$ and $E$ are chosen to be in the metric segments
$[A,C]$ and $[B,C]$, respectively, in such a~way that
$\rho(D,C)=(1-t)\rho(A,C)$ and $\rho(E,C)=(1-t)\rho(B,C)$, where
$t\in(0,1)$. Then 
\begin{equation}\label{prosta2}
\rho(D,E)\leq \frac{\sin(1-t)M}{\sin M} \; \rho(A,B).
\end{equation}
\end{lem}

The proof for the next result can be found in the proof of Theorem 4.2 in \cite{moja}.

\begin{lem}\label{lem:CAT}
  Let $x_1,x_2,x_3,x_4$ be different points of a~CAT(1) space $X$ for which $\rho(x_i,x_j)\leq M<\frac{\pi}{2}$, $i,j\in\{1,2,3,4\}$
and let $x_3\in[x_2,x_4]$ be chosen in such a~way that $d(x_3,x_4)=td(x_2,x_4)$. Then
$$
\sin^2\frac{\rho(x_1,x_3)}{2} \leq \frac{\sin(1-t)M}{\sin M}\sin^2\frac{\rho(x_1,x_4)}{2}
$$
$$
+
\frac{\sin tM}{\sin M}\frac{\max\{\cos\rho(x_2,x_4)-\cos\rho(x_2,x_1),0\}}{2}+\sin^2\frac{tM}{2}.
$$
\end{lem}

Next we present the counterpart to Lemma 2.2 of \cite{Suzuki} for CAT(1) spaces. 

\begin{lem}\label{lem:suzuki}
Let $X$ be a CAT(1) space such that $\rm{diam}(X)<\pi/2$ and let $(x_n)$ and $(y_n)$ be two sequences in $X$. Let $(\beta_n)$ 
be a~sequence in $[0,1]$ such that
\begin{equation}
 0<\liminf\beta_n\leq\limsup\beta_n<1.
\end{equation}
Suppose that $x_{n+1}=\beta_nx_n+(1-\beta_n)y_n$ for all $n\in\mathbb{N}$ and
\begin{equation}\label{granicagorna}
 \limsup \ \rho(y_{n+1},y_n)-\rho(x_{n+1},x_n)\leq 0.
\end{equation}
Then
$$
\lim \rho(x_n,y_n)=0.
$$
\end{lem}

\pf
This lemma may be proved in the same way than Lemma 2.2 in \cite{Suzuki}. This also requires to adapt Lemma 2.1 in \cite{Suzuki} to CAT(1) spaces. Both lemmas extend in an straightforward way to our situation given that the metric in a CAT(1) space of diameter smaller than $\pi/2$ is convex (see (1) of Exercise II.2.3 in \cite{bridson}), details are omitted.$\qquad\blacksquare$

\vskip.5cm

We will also need the following lemma given by Xu in \cite{oXu} on sequences of real numbers.
\begin{lem}\label{oXu}
Let $(s_n)$ be a~sequence of non-negative real numbers satisfying:
 $$s_{n+1} \leq (1 - \alpha_n)s_n + \alpha_n\beta_n + \gamma_n,\;\; n>0,$$
where $(\alpha_n)$, $(\beta_n)$ and $(\gamma_n)$ satisfy the conditions :
\begin{itemize}
\item[(i)] $(\alpha_n) \in [0, 1]$, $\sum \alpha_n = \infty$ or equivalently $\prod(1-\alpha_n)=0$,
\item[(ii)] $\limsup\beta_n \leq 0$,
\item[(iii)] $\gamma_n \geq 0$ ($n > 0$), $\sum\gamma_n < \infty$.
\end{itemize}
Then $\lim s_n = 0$.
\end{lem}

%\vspace{3mm}

To prove our main results we will need to consider the {\sl a priori bounds} given in \cite{berg} for the quadrilateral cosine on a sphere. We will not consider however the quadrilateral cosine itself but another expression which is very closely related to it and which fits better our purposes. Suppose that $X$ is a CAT($1$) space and define the~function $h$ as
$$
h(A,B;C,D)=\dfrac{\cos\rho(A,C)+\cos\rho(B,D)-\cos\rho(A,D)-\cos\rho(B,C)}{\rho(A,B)\rho(C,D)}
$$
for each four points $A,B,C,D$ of $X$ such that
$$
\max_{x,y\in\{A,B,C,D\}}d(x,y)<\pi/2 \quad \mbox{and} \quad A\neq B,\; C\neq D.
$$
Notice that $h(A,B;C,D)=h(C,D;A,B)$.

We will see that this function has some nice properties which will allow us to obtain similar bound results to the ones given by Berg and Nikolaev (compare \cite{berg}) for the quadrilateral cosine. We first give the counterpart to Lemma 2 in \cite{berg}.

\begin{lem}\label{break}
Let $A,B,C,D,X$ be chosen in such a way that $B\neq A\neq X\neq B$ and $C\neq D$. Then
$$
h(A,B;C,D)=\dfrac{\rho(A,X)}{\rho(A,B)}h(A,X;C,D)+\dfrac{\rho(X,B)}{\rho(A,B)}h(X,B;C,D).
$$
\end{lem}

\pf
The proof is straightforward.$\qquad\blacksquare$

\vskip.5cm

This result yields the next additivity property of $h$
\begin{equation}\label{sum}
h(A,B;C,D)=\dfrac{1}{nm}\sum\limits_{i=1}^n\sum\limits_{j=1}^mh(A_{i-1},A_i;C_{j-1},C_j)
\end{equation}
for $A_i\in[A,B]$ and $\rho(A_{i-1},A_i)=n^{-1}\rho(A,B)$, $i\in\{1,2,\ldots,n\}$, $A_0=A$, $A_n=B$, $C_j\in[C,D]$ and $\rho(C_{j-1},C_j)=m^{-1}\rho(C,D)$, $j\in\{1,2,\ldots,m\}$, $C_0=C$, $C_m=D$,
being the counterpart of Corollary 4 in \cite{berg}.

\vskip.5cm

The counterpart of Lemma 6 in \cite{berg} will be as follows. 

\begin{lem}\label{secondderivative}
Let distinct $P,Q,X,Y$ be given in a smooth enough CAT(1) space. Let $P_x$ be the unique point on the geodesic segment $[P,X]$ at distance $x$ from $P$ and $Q_y$ the corresponding one on the geodesic segment $[Q,Y]$ at distance $y$ from $Q$. Let $w(x,y)=\cos \rho (P_x,Q_y)$. Then, if $w(x,y)$ is in $C^2$ in a neighborhood of $(0,0)$, we have
$$
\lim\limits_{x,y\to 0^+}h(P,P_x;Q,Q_y)=\dfrac{\partial^2 w(x,y)}{\partial x\partial y}\bigg|_{(0,0)}.
$$
\end{lem}

\pf
The proof follows the same scheme of that of Lemma 6 in \cite{berg} with no significant modifications.$\qquad\blacksquare$

\vskip.5cm

Lemma 7 in \cite{berg} will turn into the following lemma for us.

\begin{lem} Let $P,Q,X,Y$ be four points in ${\mathbb S}^2$. Then, with the same notation as in Lemma \ref{secondderivative}, we have that 
\begin{equation}\label{limit}
\lim\limits_{x,y\to 0^+}h(P,P_x;Q,Q_y)=\sin\xi_x\sin\xi_y+\cos\xi_x\cos\xi_y\cos d(P,Q),
\end{equation}
where $\xi_x = \angle_{P}(Q,X)$ and $\xi_y=\pi-\angle_Q(Y,P)$.
\end{lem}

\pf
The proof of this lemma follows the same patterns than that of Lemma 7 in \cite{berg}. Consider the points $P_{x_0}$ and $Q_{y_0}$ as the respective metric projections of $P_x$ and $Q_y$ onto the line passing through $P$ and $Q$. Then, by Lemma \ref{break} applied twice, we can decompose $ h(P,P_x;Q,Q_y)$ in the following way:

\begin{align*}
d(P,P_x) & \cdot d(Q,Q_y)\cdot h(P,P_x;Q,Q_y)  \\ 
 &= d(P,P_{x_0})\cdot d(Q,Q_{y_0})\cdot h(P, P_{x_0};Q,Q_{y_0}) \\
 &+ d(P,P_{x_0})\cdot d(Q_y,Q_{y_0})\cdot h(P, P_{x_0};Q_{y_0},Q_{y}) \\
&+ d(P_x,P_{x_0})\cdot d(Q,Q_{y_0})\cdot h(P_{x_0}, P_{x};Q,Q_{y_0}) \\
&+ d(P_x,P_{x_0})\cdot d(Q_y,Q_{y_0})\cdot h(P_{x_0}, P_{x};Q_{y_0},Q_{y}) 
\end{align*}

Now we consider the same three cases than in the proof of Lemma 7 in \cite{berg} for which we obtain different results:

{\it Case (I):} Lines passing through $P$ and $P_x$, and $Q$ and $Q_y$ are both meridians and $\xi_x=\xi_y=\pi/2$. Then we obtain that 
$$\dfrac{\partial^2w(x,y)}{\partial x\partial y}\bigg|_{(0,0)}=1.$$

{\it Case (I):} If the line passing through $P$ and $P_x$ is a meridian and the one through $Q$ and $Q_y$ is equatorial ($\xi_x=\pi/2$, $\xi_y=0$), then 
$$\dfrac{\partial^2w(x,y)}{\partial x\partial y}\bigg|_{(0,0)}=0.$$

{\it Case (II):} If lines passing through $P$ and $P_x$, and $Q$ and $Q_y$ are both equatorial, i.e., $\xi_x=\xi_y=0$, then we obtain that 
$$\dfrac{\partial^2w(x,y)}{\partial x\partial y}\bigg|_{(0,0)}=w(0,0)=\cos d(P,Q).$$

Therefore, in all three cases the lemma holds. Now it is enough to follow the proof of Lemma 7 in \cite{berg} with no significant modification and where these three cases and the above decomposition of $h(P,P_x;Q,Q_y)$ are used to complete the proof of the lemma.$\qquad\blacksquare$

\vskip.5cm

From (\ref{limit}) it is clear that, in ${\mathbb S}^2$, 
\begin{equation}\label{limit2}
\lim_{B\to A, D\to C} h (A,B;C,D)\le 1
\end{equation}
whenever $B$ and $D$ move in given geodesics starting at $A$ and $C$, respectively.

Next lemma shows that it is actually possible to drop the limit condition in the above formula.

\begin{lem}\label{mainh}
Let $A,B,C,D$ be four points in ${\mathbb S}^2$ under the usual conditions, then
$$
h(A,B;C,D)\le 1.
$$
In particular, $|h(A,B;C,D)|\le 1.$
\end{lem}

\pf 
The idea for this lemma is to introduce partitions in the geodesic intervals $[A,B]$ and $[C,D]$ as small as needed. Then (\ref{sum}) and (\ref{limit2}) are considered as in Lemma 9 in \cite{berg}.$\qquad\blacksquare$

\vskip.5cm

The above considerations lead to the counterpart of the main result, Lemma 10, of Section 3.3 in \cite{berg}. 
\begin{lem}\label{maintech}
Let $X$ be a CAT(1) space. Then the inequality
$$
\left|h(A,B;C,D)\right| \leq 1
$$
holds for each foursome of points $A,B,C,D$ of $X$ such that $A\neq B$, $C\neq D$ and
$$
\max_{x,y\in\{A,B,C,D\}}\rho(x,y)<\pi/2.
$$
\end{lem}

\pf 
This result follows as a direct application of a theorem of Reshetnyak about comparison quadrangles in CAT(1) spaces. For details see the proof of Lemma 10 in \cite{berg}. $\qquad\blacksquare$

\section{Main results}

\setcounter{equation}{0}\setcounter{tw}{0}

Let us supppose that $\{b_n\}$ and $\{t_n\}$ are a~pair of sequences of nonnegative numbers satisfying the following conditions
\begin{enumerate}
  \item[(i)] $b_n,t_n\in (0,1)$ for $n\in\mathbb{N}$;
	\item[(ii)] $0<\liminf\limits_{n\to\infty} b_n \leq \limsup\limits_{n\to\infty} b_n<1$;
	\item[(iii)] $\lim\limits_{n\to\infty} t_n =0$;
	\item[(iv)] $\sum\limits_{n=1}^\infty t_n = \infty$ or equivalently $\prod\limits_{n=1}^\infty (1-t_n)=0$.
\end{enumerate}
In 2006 Chidume and Chidume (see \cite{chidume}) showed that in case of the Banach space with uniformly G\^{a}teaux differentiable norm the iterative sequence defined by 
\begin{eqnarray*}
x_1 & = & u, \\
y_{n} & = & t_n u+(1-t_n)T(x_n), \\
x_{n+1} & = & b_n x_n + (1-b_n) y_n, 
\end{eqnarray*}
where sequences $\{b_n\}$ and $\{t_n\}$ satisfy (i)--(iv),
tends to nearest with respect to $u$ fixed point of a nonexpansive mapping $T$ for which $Fix T\neq\emptyset$. Different authors have worked on this iterative scheme, some of them have modified the definition of the sequence $(y_n)$ following a Moudafi's viscosity codition (see \cite{moudafi}) as
$$
y_{n} \ = \ t_n f(x_n)+(1-t_n)T(x_n),
$$ 
where $f$ is a contraction (compare e.g. \cite{saejung}).

We present next our main result where we take up for first time the study of this viscosity iterative scheme in a nonlinear space. We particularize first for CAT(1) spaces where usual boundedness conditions on these spaces are required. Later on we will show how to adapt this result for CAT$(\kappa)$ spaces with $\kappa >0$ as well as for CAT(0) spaces where the boundedness conditions are no longer requested.

\begin{tw}
Let $C$ be a complete CAT(1) space with N--property and diameter smaller than $\pi/2$. Moreover, let $T\colon C\to C$ be a~nonexpansive mapping with $Fix T\neq\emptyset$ and let $f\colon C\to C$ be $k$-contractive with 
$$
k<\dfrac{2\sin^2\frac{M}{2}\cos M}{M^2}.$$ 
Assume further that $\rho(p,f(p)) \leq M/4$ for all $p\in Fix T$, where $M\in(0,\pi/2)$. Then there is a unique fixed point $q\in Fix T$ for which
\begin{equation}\label{projekcja}
P_{Fix T}(f(q))=q.
\end{equation}
Moreover, for each point $u\in X$ such that
$$
\rho(u,q)\leq M/4
$$ 
and for each couple of sequences $\{b_n\}$ and $\{t_n\}$ satisfying (i)--(iv), the viscosity iterative sequence defined by
\begin{eqnarray}
x_1 & = & u, \nonumber\\
y_{n} & = & t_n f(x_n)+(1-t_n)T(x_n), \label{iteracja}\\
x_{n+1} & = & b_n x_n + (1-b_n) y_n \nonumber
\end{eqnarray}
converges to $q$ satisfying \eqref{projekcja}.
\end{tw}

\pf
Let us define the function $G\colon Fix T \to Fix T$ by
$$
G(p)=P_{Fix T}(f(p)).
$$ 
From Proposition 3.5(3) in \cite{RafaAurora} we have that
\begin{equation}\label{first}
\cos \rho(G(p),f(q)) \leq \cos \rho(f(q),G(q))\cos \rho(G(p),G(q))
\end{equation}
and 
\begin{equation}\label{second}
\cos \rho(G(q),f(p)) \leq \cos \rho(f(p),G(p))\cos \rho(G(p),G(q)).
\end{equation}
Now from Lemma 3.10 we obtain
$$
\cos \rho(f(p),G(p))+\cos \rho(f(q),G(q))-\cos \rho(f((p),G(q))-\cos \rho(f(q),G(p))\leq \rho(f(p),f(q)) \cdot \rho(G(p),G(q)),
$$
so combining it with \eqref{first} and \eqref{second} yields
$$
\dfrac{\left(\cos \rho(f(p),G(p))+\cos \rho(f(q),G(q))\right)\cdot 2 \sin^2\frac{\rho(G(p),G(q))}{2}}{\rho(G(p),G(q))} \leq \rho(f(p),f(q)).
$$
Since the function $\dfrac{\sin\frac{a}{2}}{a}$ is not increasing and $f$ is a~contraction, we finally obtain
$$
\rho(G(p),G(q)) \leq \left(k\cdot\dfrac{\pi^2}{8\cos\frac{M}{4}}\right)\cdot \rho(p,q).
$$
But $k<\dfrac{\cos M}{2}$, so $G$ is also a~contraction and there is a~unique fixed point such that \eqref{projekcja} holds.

We claim next that $\displaystyle \{ x_n,y_n\colon n\in{\mathbb N}\}\subseteq B\left( q,\frac{M}{4(1-k)}\right)$. 

From the fact that each triangle of $X$ with sides smaller than $\pi/2$ is stricty convex (compare \cite[Proposition 3.1]{RafaAurora}), we get
$$
\rho (y_1,q) \leq \max\{\rho (q,T(u)),\rho (q,f(u))\}\leq \dfrac{M}{4}(1+k)<\frac{M}{4(1-k)}.
$$
By hypothesis, $\rho (q,x_1)\le M/4$. Therefore, by convexity, we also have that $\rho (q,x_2)\le \frac{M}{4(1-k)}$. For $\rho (q,y_2)$ we consider
$$
\rho(q,f(x_1))\le \rho (q,f(q))+\rho (f(q),f(x_1))\le \frac{M}{4}+k(\frac{M}{4}(1+k))\le \frac{M}{4}(1+k+k^2),
$$
and 
$$
\rho (q,T(x_1))=\rho (T(q),T(x_1))\le \rho (q,x_1)\le \frac{M}{4}(1+k+k^2).
$$
So, by convexity again, $\rho (q,y_2)\le \frac{M}{4}(1+k+k^2)$. Continuing in this way our claim follows. In the process we have also shown that the same holds for sequences $\{f(x_n)\}$ and $\{ T(x_n)\}$ so, in particular, it also follows from the convexity of the metric that
\begin{equation}\label{convexbound}
d(x,y) \leq \rho (x,q)+\rho (q,y) \leq 2\dfrac{M}{4(1-k)} \leq M,
\end{equation}
for each pair of points $x,y \in \overline{conv}\{x_1,x_2,\ldots,T(x_1),T(x_2),\ldots,f(x_1),f(x_2),\ldots\}$ (the definition and some properties of this convex hull may be found in \cite[Proposition 4.5]{RafaAurora}).

{\bf Step 1.} We prove first that $\lim_{n\to \infty}\rho(x_n,T(x_n))=0$. By applying Lemma \ref{CAT} twice for triangles $\Delta(f(x_n),T(x_n),T(x_{n+1}))$ and $\Delta(f(x_n),f(x_{n+1}),T(x_{n+1}))$, respectively, as well as assumptions on $T$ and $f$, it follows that
$$
\rho(y_n,y_{n+1})\leq \frac{\sin(1-t_n)M}{\sin M}\rho(T(x_n),T(x_{n+1}))+\frac{\sin t_nM}{\sin M}\rho(f(x_n),f(x_{n+1}))+|t_{n}-t_{n+1}|M
$$
$$
\leq \left(\frac{\sin(1-t_n)M}{\sin M}+k\frac{\sin t_nM}{\sin M}\right)\rho(x_n,x_{n+1})+|t_{n}-t_{n+1}|M,
$$
what leads to
%$$
%\rho(y_n,y_{n+1})-\rho(x_n,x_{n+1}) \leq \dfrac{2kM\sin\frac{t_nM}{2}}{\sin M}+|t_{n}-t_{n+1}|M
%$$
$$
\rho(y_n,y_{n+1})-\rho(x_n,x_{n+1}) \leq  \left(\frac{\sin(1-t_n)M}{\sin M}+k\frac{\sin t_nM}{\sin M}-1\right)\rho(x_n,x_{n+1})+|t_{n}-t_{n+1}|M.
$$
Taking limit when $n$ goes to infinity makes the right hand side convergent to $0$ which yields \eqref{granicagorna}. From Lemma \ref{lem:suzuki} we have now that $\rho(y_n,x_n) \to 0.$ At the same time,
$$
\rho(y_n,T(x_n)) = \rho(t_nf(x_n)+(1-t_n)T(x_n),T(x_n)) = t_n \rho(f(x_n),T(x_n)) \leq t_n M,
$$
what leads to
$$
\rho(y_n,T(x_n)) \to 0.
$$
Finally, also
$$
\rho(x_n,T(x_n)) \leq \rho(x_n,y_n) + \rho(y_n,T(x_n))
$$
tends to 0 if $n\to\infty$.

{\bf Step 2.} In this step we show how to apply Lemma \ref{oXu} to our scheme. Let us consider a~triangle $\Delta(q,x_n,y_n)$. Since $x_{n+1}\in[x_n,y_n]$, from Lemma \ref{lem:pierwszy}, it follows that
$$
\sin^2\frac{\rho(q,x_{n+1})}{2}\leq b_n\sin^2\frac{\rho(q,x_n)}{2}+(1-b_n)\sin^2\frac{\rho(q,y_n)}{2}=(*).
$$
Now let us take a~triangle $\Delta(q,f(x_n),T(x_n))$ and apply Lemma \ref{lem:CAT} to get
$$
(*)\leq b_n\sin^2\frac{\rho(q,x_n)}{2}+(1-b_n)\left(\frac{\sin(1-t_{n})M}{\sin M}\sin^2\frac{\rho(T(x_n),q)}{2}+\right.
$$ 
$$
\left.
+\frac{\sin t_{n}M}{\sin M}\frac{\max\{\cos \rho(f(x_n),T(x_n))-\cos \rho(f(x_n),q),0\}}{2}+\sin^2\frac{t_nM}{2}\right)
$$
$$
\le b_n\sin^2\frac{\rho(q,x_n)}{2}+(1-b_n)\left(\frac{\sin(1-t_{n})M}{\sin M}\sin^2\frac{\rho(T(x_n),q)}{2}+\right.
$$
$$
+\frac{\sin t_{n}M}{\sin M}\frac{\max\{\cos \rho(f(x_n),T(x_n))-\cos \rho(f(x_n),q)+\cos \rho(f(q),q)-\cos \rho(f(q),T(x_n)),0\}}{2}+
$$
$$
\left.
+\frac{\sin t_{n}M}{\sin M}\frac{\max\{\cos \rho(f(q),T(x_n))-\cos \rho(f(q),q),0\}}{2}+\sin^2\frac{t_nM}{2}\right)
$$
by Lemma \ref{maintech},
$$
\leq b_n\sin^2\frac{\rho(q,x_n)}{2}+(1-b_n)\left(\frac{\sin(1-t_{n})M}{\sin M}\sin^2\frac{\rho(T(x_n),q)}{2}+\right.
$$
$$
+\frac{\sin t_{n}M}{\sin M}\frac{\rho(f(x_n)f(q))\rho(T(x_n),q)}{2}
$$
$$
\left.
+\frac{\sin t_{n}M}{\sin M}\frac{\max\{\cos \rho(f(q),T(x_n))-\cos \rho(f(q),q),0\}}{2}+\sin^2\frac{t_nM}{2}\right).
$$
Since $T$ is a nonexpansive mapping, $q$ is a fixed point of $T$, $f$ is a contraction and the mapping $x\mapsto \frac{x}{\sin\frac{x}{2}}$ is increasing, 
$$
\sin^2\frac{\rho(q,x_{n+1})}{2}\leq
\sin^2\frac{\rho(q,x_n)}{2}\left(b_n+(1-b_n)\frac{\sin(1-t_{n})M}{\sin M}+(1-b_n)k\frac{\sin t_{n}M}{\sin M}\frac{M^2}{2\sin^2\frac{M}{2}}\right)+
$$
$$
+(1-b_n)\left(
\frac{\sin t_{n}M}{\sin M}\frac{\max\{\cos \rho(f(q),T(x_n))-\cos \rho(f(q),q),0\}}{2}+\sin^2\frac{t_nM}{2}\right).
$$

Fix $\theta\in (0,1)$ such that $\displaystyle k=\theta \cos M\frac{2\sin^2\frac{M}{2}}{M^2}$. Then
$$
\begin{array}{lll}
\dfrac{\sin(1-t_n)M}{\sin M}+k\dfrac{\sin t_nM}{\sin M}\dfrac{M^2}{2\sin^2\frac{M}{2}} & = &
\dfrac{\sin(1-t_n)M}{\sin M}+\dfrac{\sin t_nM}{\sin M}\theta\cos M
\\
& = & \dfrac{\cos t_nM\sin M - (1-\theta)\sin t_nM\cos M}{\sin M} \ = \ (*).  
\end{array}
$$
Recalling the concavity of sine function, $t_n\sin M <\sin t_nM$, and so
$$
(*) \le 1-t_n (1-\theta)\cos M.
$$

Let us choose now a positive number $A$ in such a way that 
$$
A<\cos M - k\frac{M^2}{2\sin^2\frac{M}{2}},
$$
hence $A< (1-\theta )\cos M$, and fixing $b\in (0,1)$ such that  $b_n\le b$ for all $n$ large enough
$$
b_n+(1-b_n)\left(\dfrac{\sin(1-t_n)M}{\sin M}+k\dfrac{\sin t_nM}{\sin M}\dfrac{M^2}{2\sin^2\frac{M}{2}}\right)
\leq 1 - (1-b_n)A t_n \leq 1 - (1-b)A t_n.
$$
So we obtain that
$$
\sin^2\frac{\rho(q,x_{n+1})}{2}\leq
\sin^2\frac{\rho(q,x_n)}{2}(1-(1-b)A t_n)+
$$
$$
+t_n M(1-b_n)\left(\frac{1}{\sin M}\frac{\max\{\cos \rho(f(q),T(x_n))-\cos \rho(f(q),q),0\}}{2}+\frac{1}{2}\sin\frac{t_nM}{2}\right).
$$
Taking $s_n=\sin^2\dfrac{\rho(q,x_n)}{2}$, $\alpha_n=(1-b)A t_n$, on account of Lemma \ref{oXu}, it suffices to show that
$$
\beta_n=\dfrac{M(1-b_n)}{A(1-b)}\left(\frac{1}{\sin M}\frac{\max\{\cos \rho(f(q),T(x_n))-\cos \rho(f(q),q),0\}}{2}+\frac{1}{2}\sin\frac{t_nM}{2}\right)
$$
tends to $0$, that is, that
\begin{equation}\label{step3}
\limsup_{n\to\infty}\cos \rho(f(q),T(x_n)) \leq \cos \rho(f(q),q).
\end{equation}

{\bf Step 3.} Now we will show (\ref{step3}).
Let us consider the~set of projections $\{P_{[f(q),q]}(y_n) \colon$ $n\in\mathbb{N}\}$. Suppose that there is a~subsequence $(y_{k_n})$ of $(y_n)$ such that
\begin{equation}\label{projekcje}
P_{[f(q),q]}(y_{k_n})\in [p_1,p_2],
\end{equation}
where $p_1\neq q\neq p_2$. This sequence has a~regular subsequence with respect to $\Delta$--convergence denoted again by $(y_{k_n})$ (see Corollary 4.4 in \cite{RafaAurora} and \cite{delta}).
Now from Proposition 4.5 in \cite{RafaAurora} we get that there is a~unique asymptotic center $A:=A_{C}(y_{k_n})$ which belongs to $\bigcap_{n=1}^\infty \overline{conv}\{y_{k_n},y_{k_{n+1}},\ldots\}$. 
Combining this with \eqref{projekcje} and the N--property of the space we obtain that 
\begin{equation}\label{projekcjaa}
 P:=P_{[f(q),q]}(A)\in[p_1,p_2].
\end{equation}
On the other hand from Proposition 3.13 in \cite{wspolna} it follows that $A$ is a~fixed point of $T$. 

Now let us consider a~comparison triangle $\Delta(\bar{A},\bar{q},\bar{P})$ of $\Delta(A,q,P))$. Let us suppose that an~angle $\angle_{\bar{q}}(\bar{P},\bar{A})$ is not smaller than $\pi/2$. Then we get that in the comparison triangle
$$
\cos d(\bar{A},\bar{P})= \cos d(\bar{A},\bar{q})\cos d(\bar{P},\bar{q})+\sin d(\bar{A},\bar{q})\sin d(\bar{P},\bar{q})\cos\angle_{\bar{q}}(\bar{P},\bar{A})
$$
$$
\leq \cos d(\bar{A},\bar{q})\cos d(\bar{P},\bar{q}).
$$
But $\rho(A,P)=d(\bar{A},\bar{P})$ and $\rho(A,q)=d(\bar{A},\bar{q})$. Therefore
$$
\cos\rho(A,P)\leq \cos\rho(A,q)\cos d(\bar{P},\bar{q})< \cos\rho(A,q),
$$
what contradicts the definition of $P$ in \eqref{projekcjaa}. Comparing angles in the triangle $\Delta(A,q,P)$ and its comparison triangle we obtain that 
$$
\angle_q(P,A) \leq \angle_{\bar{q}}(\bar{P},\bar{A})<\frac{\pi}{2}
$$
(see \cite[Proposition 1.7 (4), p. 161]{bridson}).

Now it suffices to notice that from the Alexandrov definition of angles it follows that the angle $\angle_q(A,P)$ in the triangle $\Delta(A,q,P)$ and $\angle_q(A,f(q))$ in $\Delta(A,q,f(q))$ are equal, what yields to
$$
\angle_q(A,f(q))<\frac{\pi}{2}.
$$
On the other hand $q$ is the projection of $f(q)$ onto the closed $D_1$--convex set $Fix(T)$ of fixed points of $T$,
what on account of condition (2) in Proposition 3.5 in \cite{RafaAurora} follows to the estimation
$$
\angle_q(A,f(q))\geq \frac{\pi}{2},
$$
a contradiction.

Consequently our guess on the existence of $p_1$ and $p_2$ is wrong and so  
$$
\liminf \rho(y_{n},f(q)) \geq \liminf \rho(P_{[q,f(q)]}y_{n},P_{[q,f(q)]}f(q))= \lim \rho(P_{[q,f(q)]}y_{n},f(q))= \rho(q,f(q)).
$$
Since $t_n\to 0$, the following inequality
$$
\liminf \rho(T(x_n),f(q)) \geq \rho(q,f(q))
$$
also holds, what completes the proof.$\qquad\blacksquare$

\vspace{3mm}

Since the distances in model spaces $M_\kappa^2$ with $\kappa>0$ are obtained by scaling the distances on the unit sphere, we can extend the previous result to general CAT$(\kappa)$ spaces with $\kappa>0$ as the following theorem shows.

\begin{tw}
Let $C$ be a complete CAT($\kappa$) space, $\kappa>0$, with N--property and diameter smaller than $D_{\kappa}/2$. Moreover, let $T\colon C\to C$ is a nonexpansive mappings with $Fix T\neq\emptyset$ and $f\colon C\to C$ a $k$-contraction $k<\dfrac{2\sin^2\frac{M}{2}\cos M}{M^2}$ and such that $\rho(p,f(p)) \leq M/(4\sqrt{\kappa})$, $p\in Fix T$, for a fixed number $M\in(0,\pi/2)$. Then there is the unique fixed point $q\in Fix T$ for which \eqref{projekcja} holds.

Moreover, for each point $u\in X$ such that
$$
d(u,q)\leq M/(4\sqrt{\kappa})
$$ 
and for each couple of sequences $\{b_n\}$ and $\{t_n\}$ satisfying (i)--(iv), the viscosity iterative sequence defined by \eqref{iteracja} converges to $q$.
\end{tw}

\pf
Using the model space $M^2_\kappa$ defined in Section 2 we get
$$
d(x,y) = \sqrt{\kappa}
d_{M^2_1}(x,y).
$$
Now to obtain our result it sufices to repeat all steps from the proof of Theorem 4.1.$\qquad\blacksquare$

\vspace{3mm}

Let us notice that each CAT(0) space is also a CAT($\kappa$) space for all $\kappa>0$, so the above results may be modified to this class of spaces.

\begin{tw}
Let $C$ be a complete CAT(0) space with N--property. Let $T\colon C\to C$ be a nonexpansive mapping with $Fix T\neq\emptyset$ and let $f\colon C\to C$ be $k$-contraction $k<\dfrac{1}{2}$. Then there is the unique fixed point $q\in Fix T$ for which \eqref{projekcja} holds.

Moreover, for each point $u\in X$ 
and for each couple of sequences $\{b_n\}$ and $\{t_n\}$ satisfying (i)--(iv), the viscosity iterative sequence defined by \eqref{iteracja} converges to $q$.
\end{tw}

\pf
Since each CAT(0) space is also a CAT($\kappa$) space for all positive $\kappa$ (see \cite[Theorem 1.12,
p. 165]{bridson}), given $u$ and $T$ as in the statement we can choose $\kappa$ and $M$ satisfying assumptions of Theorem 4.2.
Now
we just need to apply that theorem for such $\kappa$ and the result follows.$\qquad\blacksquare$

\vspace{3mm}

\section{Appendix: A CAT(0) space failing the N-property}
\setcounter{equation}{0}

In our results we have required the space to have property N. This property has been recently introduced in \cite{RafaAurora} where it was claimed that property N was to be very common within the class of CAT$(\kappa)$ spaces but no such space failing property N was provided. We show next a CAT$(0)$ space failing this property. As long as the author knows, this is the first example of that class in the literature.

\vspace{3mm}

\noindent{\bf Example 5.1}
Let us consider two flat triangles in $\mathbb{E}^3$
$$
A=(-1,0,4); B=(1,0,4); C=(0,0,0)
$$
and
$$
C=(0,0,0); D=(0,0,4); E=\left(0,\frac{3\sqrt{7}}{8},\frac{1}{8}\right).
$$

Now let us define a~distance in the space $X=\Delta(A,B,C)\cup\Delta(C,D,E)$ as the lenght of the shortest path in $X$ connecting two points.
To check that our space is CAT$(0)$ it suffices to notice that $X$ is a gluing space of two subsets of $\mathbb{E}^2$, so one may apply Reshetnyak's Gluing Theorem (see e.g. \cite[Theorem 9.1.21, p. 316]{burago}).

We will prove that 
\begin{equation}\label{wzorek}
P_{[C,E]}(D)=P_{[C,E]}\left(\frac{1}{2}A+\frac{1}{2}B\right)\not\in\left[P_{[C,E]}(A),P_{[C,E]}(B)\right].
\end{equation}

Clearly, the projection of $D$ is equal to $F=\dfrac{1}{2}C+\dfrac{1}{2}E=\left(0,\dfrac{3\sqrt{7}}{16},\dfrac{1}{16}\right)$ because $\angle_F(D,E)=\dfrac{\pi}{2}$.
Now let us calculate
$$
d(A,C)=\sqrt{1^2+4^4}=\sqrt{\frac{136}{8}}
$$
and
$$
d(A,F)=\sqrt{\left(\frac{3\sqrt{7}}{16}+1\right)^2+\left(\frac{1}{16}-4\right)^2}=\sqrt{\frac{134+3\sqrt{7}}{8}}>\sqrt{\frac{136}{8}}.
$$
Because of symmetry $P_{[C,E]}(A)=P_{[C,E]}(B)\neq F$. This completes the proof of \eqref{wzorek}.$\qquad\blacksquare$

\vspace{3mm}

We close this work by raising the following open problem.

{\sl Open problem: } Can we drop the condition on the N-property from our main results?

\end{document}